

\magnification=1200
\nopagenumbers
\parindent= 15pt
\baselineskip=14pt

\hsize=13.65cm
\vsize=19cm
\hoffset=-0.15cm
\voffset=0.3cm

\input amssym.def
\input amssym.tex

 at6.5pt
\font\srm=cmr8 

\font\cmx=cmbxti10
\font\csc=cmcsc10
\font\title=cmr12 at 14pt

\font\teneusm=eusm10    
\font\seveneusm=eusm7  
\font\fiveeusm=eusm5    
\newfam\eusmfam
\def\eusm{\fam\eusmfam\teneusm}
\textfont\eusmfam=\teneusm 
\scriptfont\eusmfam=\seveneusm
\scriptscriptfont\eusmfam=\fiveeusm

\def\Re{{\rm Re}\,}
\def\Im{{\rm Im}\,}
\def\sgn{{\rm sgn}\,}
\def\txt#1{{\textstyle{#1}}}
\def\scr#1{{\scriptstyle{#1}}}

\def\r#1{{\rm #1}}
\def\e#1{{\eusm #1}}
\def\varGamma{{\mit\Gamma}}
\def\B#1{{\Bbb #1}}

\def\sgn{{\rm sgn}}

\def\rightheadline{\hfil{\srm On
sums of Hecke--Maass eigenvalues
squared over primes in short intervals}
\hfil\tenrm\folio}
\def\leftheadline{\tenrm\folio\hfil{\srm Y. Motohashi}\hfil}
\def\emptyheadline{}
\headline{\ifnum\pageno=1 \emptyheadline\else
\ifodd\pageno \rightheadline \else \leftheadline\fi\fi}

\def\firstpage{\hss{\vbox to 1cm
{\vfil\hbox{\rm\folio}}}\hss}
\def\emptyfootline{\hfil}
\footline{\ifnum\pageno=1\firstpage\else
\emptyfootline\fi}

\centerline{\title On Sums of
Hecke--Maass Eigenvalues Squared}
\medskip
\centerline{\title over Primes in Short
Intervals} 
\vskip 0.7cm
\centerline{\csc Yoichi Motohashi}
\vskip 1cm 
\noindent
{\bf Abstract:}  We prove a uniform estimate
for sums of Hecke--Maass eigenvalues squared over 
primes in short intervals; the precise assertion is
given in Theorem 1 that can be regarded as an analogue,
for {\it all\/} Maass forms,
of Hoheisel's theorem [5] on the
existence of primes in short intervals. Our argument is modelled
after our treatment [11] of Linnik's least prime number
theorem for arithmetic progressions.
We stress that constants in the present work, 
including those implicit, are all {\it universal and 
effectively computable\/}; we shall neither mention this repeatedly
nor pay any particular attention to numerical precision of 
our constants.
\bigskip
\noindent
{\bf 1. Introduction}  
\smallskip
\noindent
We begin with a brief about our normalisation of
automorphic forms: the article [14] contains complete
details. Thus,
we work on the Lie group $\r{G}=\r{PSL}(2,\B{R})$ 
equipped with the 
co-ordinate system $\r{G}=\r{NAK}$,
where $\r{N}=\Big\{\Big[{1\atop}{x\atop1}\Big]: x\in
\B{R}\Big\}$, $\r{A}=\Big\{\Big[{\sqrt{y}\atop}
{\atop1/\sqrt{y}}\Big]: y>0\Big\}$ 
and $\r{K}=\Big\{\Big[{\hfill\cos\xi\atop-\sin\xi}
{\sin\xi\atop\cos\xi}\Big]: \xi\in\B{R}/\pi\B{Z}\Big\}$. 
The Casimir operator is defined in the form 
$\Omega=-y^2\big(\partial_x^2+
\partial^2_y\big)+y\partial_x\partial_\xi$. We let
$L^2(\varGamma\backslash\r{G})$ with $\varGamma=
\r{PSL}(2,\B{Z})$ stand for 
the Hilbert space composed of 
all left $\varGamma$-automorphic functions on
$\r{G}$ which are square integrable over any fundamental
domain of $\varGamma$ against the invariant measure
$dxdyd\xi/\pi y^2$. This is a unitary representation
space of $\r{G}$, since all right translations by
the elements of $\r{G}$ induce unitary maps there. 
In this construction, let $V$ be an irreducible representation or
subspace occurring in $L^2(\varGamma\backslash\r{G})$.
We assume that $V$ be cuspidal so that
the Casimir operator restricted to $V$ is a constant multiplication:
$\Omega|_V=\big({1\over4}-\nu_V^2\big)\cdot1$; the
constant $\nu_V$ is often termed the spectral data of $V$.
Throughout the present work, it is assumed that
$$
\hbox{$V$ belong
to the unitary principal series; hence, $\nu_V\in i\B{R}$.}
\eqno(1.1)
$$
Then, $V$ is generated, via applications of left invariant
differential operators, by a Maass wave:
$$
\psi_V(x+iy)={2\pi^{\nu_V+1/2}\over
\Gamma\big(\nu_V+{1\over2}\big)}
\sum_{n\ne0}\varrho_V(n)y^{1/2}K_{\nu_V}(2\pi|n|y)
\exp(2\pi inx),\eqno(1.2)
$$
where $x+iy$ corresponds to a generic
point on $\r{G}/\r{K}$, the hyperbolic upper half-plane, and
$K_\nu$ is the $K$-Bessel function of order $\nu$. What is
essential in this expansion is that the sequence
$\{\varrho_V(n): \B{Z}\ni n\ne0\}$ depends only
on $V$, save for constant multipliers of unit absolute value.
We impose further that $V$ be 
Hecke invariant and $\psi_V$ be
of unit length as a vector of
$L^2(\varGamma\backslash\r{G})$. More precisely,
the Hecke operator $\r{T}(n)$ for
each $n\in\B{N}$  is defined to be the map $f(\r{g})\mapsto
n^{-1/2}\sum_{ad=n}\sum_{b\bmod d}f\Big(\left[{a\atop}
{b\atop d}\right]\r{g}\Big)$, with $a,b,d$ being
non-negative integers; and 
the $n$-th Hecke--Maass eigenvalue is denoted by
$\tau_V(n)$: $\r{T}(n)|_V=\tau_V(n)\cdot1$. In addition,
one may assume that $V$ be invariant with respect to
the involution $(x,y,\xi)\mapsto(-x,y,-\xi)$, which
induces the parity sign $\epsilon_V=\pm1$. Summing up,
we have the following
normalisation of the Fourier coefficients of $\psi_V$ or rather
of the representation $V$:
$$
\varrho_V(n)=\epsilon_V^{(1-\sgn(n))/2}\varrho_V(1)
\tau_V(|n|).\eqno(1.3)
$$
As is well-known, Kim--Sarnak [6, Appendix 2]
proved the bound 
$|\tau_V(n)|\le d(n)n^{7/64}$, where $d(n)$ is
the number of divisors of $n$. 
\medskip
With this, the Hecke $L$-function associated with $V$ is defined
to be
$$
\eqalignno{
L(s;V)&\,=\sum_{n=1}^\infty\tau_V(n)n^{-s}\cr
&\,=\prod_p\bigg[\Big(1-\alpha_V(p)p^{-s}\Big)
\Big(1-\alpha^{-1}_V(p)p^{-s}\Big)\bigg]^{-1},&(1.4)
}
$$
where $p$ stands for a generic prime and 
$$
\hbox{either $|\alpha_V(p)|=1$ or $1<\pm\alpha_V(p)
\le p^{7/64}$}.\eqno(1.5)
$$
The sum and the product 
converge absolutely for $\Re s>1$, and $L(s;V)$ continues to
an entire function, satisfying the functional equation
of the Riemann type
$$
\Gamma_2(s)L(s;V)=\kappa_{V,2}
\overline{\Gamma_2(1-\bar{s})}L(1-s;V),\eqno(1.6)
$$
where $|\kappa_{V,2}|=1$ and
$$
\Gamma_r(s)=\pi^{-rs/2}
\prod_{j=1}^r\Gamma\big(\txt{1\over2}
(s+\gamma^{(r,j)}_V)\big),\quad 
\sum_{j=1}^r|\gamma_V^{(r,j)}|\ll|\nu_V|.
\eqno(1.7)
$$
We use an abbreviated notation for
$\Gamma$-factors. The dependency on specific 
representations is made implicit: essential for our
purpose is the size of the shift parameters $\gamma_V^{(r,j)}$.
\par
The Rankin $L$-function associated with $V$ is defined to be
$$
L(s;V\times V)=\zeta(2s)
\sum_{n=1}^\infty\tau_V^2(n)n^{-s}.\eqno(1.8)
$$
This converges absolutely for $\Re s>1$ and
continues to a meromorphic function which
is regular except for the
simple pole at $s=1$ with the residue 
${1\over2}|\varrho_V(1)|^{-2}$ and satisfies
the functional equation
$$
\Gamma_4(s)L(s;V\times V)=\kappa_{V,4}
\overline{\Gamma_4(1-\bar{s})}
L(1-s;V\times V),\eqno(1.9)
$$
with $|\kappa_{V,4}|=1$. 
\par
The finiteness of the order 
of these $L$-functions is
well-known; see, e.g., [13, Vol.\ II, Section 5.1] 
for a precise treatment of the issue.
As a consequence, follow convexity polynomial
bounds for these $L$-functions which are uniform in all
involved parameters in any fixed vertical strip of the
$s$-plane; a further explanation 
is to be given in the later part of the next section.
\medskip
With this, our aim is to establish
\medskip
\noindent
{\bf Theorem 1.} {\it  Under the convention 
$(1.1)$--$(1.3)$,
there exist constants $c_0, \theta_0>0$ such that 
we have, uniformly
for $(\log x)^{-1/2}\le\theta\le\theta_0$,
$|\nu_V|^{1/\theta}\le x$,
$$
\sum_{x-y\le p\le x}\tau_V^2(p)=  {y\over \log x}
\Big(1+O\big(e^{-c_0/\theta}\big)\Big),
\quad y=x^{1-\theta}.\eqno(1.10)
$$
In particular, if $\theta$ is
small, then the sum on the left
is proportional to the number of primes in
the interval $[x-y,x]$ for any $V$.
}
\medskip
\noindent
The novelty of this assertion 
is in that $c_0,\theta_0$ and the implied constants are
all universal and effectively computable.
In (i) of Section 5 is a comparison between $(1.10)$
and Hoheisel's prime number theorem for short intervals.
Our theorem is a somewhat over-simplified 
version of what the argument of the present work is
actually capable to yield. 
Especially, the lower bound for $x$ could be made more flexible. 
However, we will not give details, since they should
transpire readily from our proof of the theorem in Section 4
and since the saliency of our assertion
lies in the very fact that an estimate of
the above type is attainable with $x$ which is of an
effective polynomial order in $|\nu_V|$.
\medskip
As is well-known, the fact that 
$L(s; V\times V)$ does not have exceptional zeros
or more precisely that the function does not have any
real zero within an effective 
distance $O\big((\log|\nu_V|)^{-1}\big)$
from $s=1$
was discovered by Goldfeld--Hoffstein--Lieman--Lockhart
[4, Appendix]; and the associated uniform zero-free region
$$
L(s;V\times V)\ne0,\quad \Re s>
1-{a\over\log(|\nu_V|+|\Im s|))}\,,\eqno(1.11)
$$
with a constant $a>0$,  was proved by
Ramakrishnan--Wang [17]: the assertion $(2.12)$ below
contains this. Accordingly, what we really need 
in order to prove $(1.10)$ is the following
zero-density estimate of the Linnik type:
\medskip
\noindent
{\bf Theorem 2.} {\it Let $N_{V\times V}(\alpha, T)$ with 
$T\ge1$ be the number of zeros of $L(s;V\times V)$ 
in the rectangular region
${1\over2}\le\Re s\le\alpha\le 1$, $|\Im s|\le T$. 
Then, there exists a constant $\omega>0$ such that
$$
N_{V\times V}(\alpha, T)
\ll(|\nu_V|T)^{\omega(1-\alpha)}.\eqno(1.12)
$$
}
\par
\noindent
We shall develop a proof of this in the next three sections.
\medskip
\noindent
{\csc Acknowledgements:} We are grateful to
S. Gelbart, J. Hoffstein, A. Ivi\'c, M. Jutila, A. Perelli,
D. Ramakrishnan, A. Sankaranarayanan and N. Watt for their
kind comments and help. 
\bigskip
\noindent
{\bf 2. Symmetric power {\cmx L}-functions} 
\par
\noindent
Our argument to prove $(1.12)$ is similar, to
an extent, to that of our recent work [15]: we rely on the
theory of symmetric power $L$-functions and the 
$\Lambda^2$-sieve method.
\smallskip
Deferring the clarification of the convergence and
the regularity of relevant $L$-functions, we start with the 
Shimura factorisation
$$
L(s;V\times V)=\zeta(s)L(s;\r{sym}^2V),\eqno(2.1)
$$
where
$$
\eqalignno{
L(s;\r{sym}^\ell V)
=&\,\prod_p\prod_{j=0}^\ell
\bigg(1-\alpha^{\ell-2j}_V(p)p^{-s}\bigg)^{-1}\cr
=&\,\prod_p\exp\Bigg(\sum_{h=1}^\infty
\sum_{j=0}^\ell\alpha^{h(\ell-2j)}_V(p)/
hp^{hs}\Bigg)\cr
=&\,\sum_{n=1}^\infty\tau_V^{(\ell)}(n)n^{-s}
,&(2.2)
}
$$
say; the definition $(1.4)$ is the same as the convention 
$V=\r{sym}^1V$ and $\tau_V(n)=\tau_V^{(1)}(n)$.
We note that $\tau_V^{(\ell)}(n)$ is real, since
$\sum_{j=0}^\ell\alpha_V^{h(\ell-2j)}(p)$ is so; see $(1.5)$.
The identity $(2.1)$
means that because of
the well-known theory on zeros of $\zeta(s)$
it suffices for us to prove instead
$$
N_{\r{sym}^2V}(\alpha, T)\ll 
(|\nu_V|T)^{\omega(1-\alpha)},\eqno(2.3)
$$
with the left side being an obvious
analogue of $N_{V\times V}(\alpha,T)$ 
for $L(s;\r{sym}^2V)$. In the light of the discussion
developed in [15, Part I], we are led to
the $\Lambda^2$-sieve situation
$$
\sum_{n\le N}\big(\tau^{(2)}_V(n)\big)^2
\bigg(\sum_{d|n}
\lambda_d\bigg)^2,\eqno(2.4)
$$
where the real numbers $\{\lambda_d\}$ are supported
on the set of square-free integers and such that
$\lambda_1=1$, and $\lambda_d=0$ for $d>R$ with
a large $R$. Then, we note that
$$
\big(\tau^{(\ell)}_V(n)\big)^2\le
\tau_{V\times V}^{(\ell)}(n),\eqno(2.5)
$$
in which the right side is the coefficient of the Rankin
convolution
$$
\eqalignno{
L(s;\r{sym}^\ell V\times\r{sym}^\ell V)
=&\,\prod_p\prod_{j=0}^\ell\prod_{k=0}^\ell
\Big(1-\alpha_V^{2(\ell-j-k)}(p)p^{-s}\Big)^{-1}\cr
=&\,\prod_p\exp\Bigg(\sum_{h=1}^\infty
\bigg(\sum_{j=0}^\ell\alpha^{h(\ell-2j)}_V(p)
\bigg)^2/hp^{hs}\Bigg)\cr
=&\,\sum_{n=1}^\infty\tau_{V\times V}^{(\ell)}(n)
n^{-s}.&(2.6)
}
$$
In particular, $\tau_{V\times V}^{(\ell)}(n)$ are all non-negative;
this trivial observation will play an important r\^ole in the proof
of Lemmas 1 and 2 below.
To prove the assertion $(2.5)$,
we first expand the exponentiated
$p$-term in the middle line of $(2.2)$ and get
$$
\tau_V^{(\ell)}(p^m)=\sum_{k=0}^m{1\over k!}
\sum_{h_1+\cdots+h_k=m}\,
\prod_{r=1}^k\bigg({1\over h_r}\sum_{j=0}^\ell
\alpha^{h_r(\ell-2j)}_V(p)\bigg),\eqno(2.7)
$$
with an obvious restriction on variables. Since 
$\tau_{V\times V}^{(\ell)}(p^m)$ has the 
same construction as this,
excepting that the sums over $j$ are all squared, the Cauchy inequality
gives
$$
\big(\tau_V^{(\ell)}(p^m)\big)^2\le\tau_{V\times V}^{(\ell)}(p^m)
\sum_{k=0}^m{1\over k!}
\sum_{h_1+\cdots+h_k=m}
{1\over h_1\cdots h_k}.\eqno(2.8)
$$
The sum on the right is the coefficient of $p^{-ms}$
in the expansion of
$\exp\big(-\log(1-p^{-s})\big)$, which confirms $(2.5)$.
\par
Hence, we consider, instead of $(2.4)$,
$$
\sum_{n\le N}\tau^{(2)}_{V\times V}(n)
\bigg(\sum_{d|n}\lambda_d\bigg)^2;\eqno(2.9)
$$
the merit of the replacement is to be felt at 
$(3.11)$ and $(3.23)$ below. In this context, 
we exploit the factorisation
$$
L(s;\r{sym}^\ell V\times\r{sym}^\ell V)=
\zeta(s)\prod_{l=1}^\ell L(s;\r{sym}^{2l}V)\eqno(2.10)
$$
which can be shown by the second lines of $(2.2)$ 
and $(2.6)$ and by the identity
$$
\sum_{j=0}^\ell\sum_{k=0}^\ell X^{2(\ell-j-k)}
=\bigg({X^{\ell+1}-X^{-(\ell+1)}\over X-X^{-1}}\bigg)^2
=\sum_{l=0}^\ell\sum_{m=0}^{2l} X^{2(l-m)}.\eqno(2.11)
$$
Our discussion in the sequel is based on the following:
\medskip
\noindent
{\bf Lemma 1.} {\it The functions $L(s;\r{sym}^{2l}V)$,
$l=1,\,2$, are entire and satisfy functional
equations with the $\Gamma$-factors
$\Gamma_3(s)$ and $\Gamma_5(s)$, 
respectively. Also
$$
L(s;\,\r{sym}^2 V\times\r{sym}^2 V)\ne0,
\quad \Re s>1-{a\over\log(|\nu_V|+|\Im s|)}\,,
\eqno(2.12)
$$
with the same constant $a>0$ as in $(1.11)$. In particular, it
holds that
$$
{L'\over L}(1;\r{sym}^{2l} V)\ll \log|\nu_V|,
\quad l=1,\,2.
\eqno(2.13)
$$
\/}
\par
\noindent
{\it Proof\/}. The first assertion stems from
Gelbart--Jacquet [2] and Kim--Shahidi [7], respectively;
in fact, we ought to cite 
relevant works more but we restrict ourselves to 
those the most directly related to our present purpose.
The uniform zero-free region $(2.12)$ is due to 
Ramakrishnan--Wang 
[17, Theorem 4.12]. The bound $(2.13)$ can
be deduced via $(2.12)$ or $L(s;\r{sym}^{2l}V)\ne0$, 
$l=1,2$, in the region indicated there by a 
well-known argument of Landau [9] which is
a fine application of the Borel--Carath\'eodory convexity theorem;
see also [13, Vol.\ I, Lemma 1.4] and
[19, Section 3.9]. We shall explain salient points
of its adaptation to the present situation. Thus, we shall first show that
for $\ell=1,2$
$$
L(s;\r{sym}^\ell V\times \r{sym}^\ell V)^{-1}
\ll((\sigma-1)^{-1}+1)
|\nu_V|^c,\quad \sigma=\Re s>1;\eqno(2.14)
$$
incidentally, we stress that here and
in what follows the symbol $c$ is to stand for a constant in
the sense of Abstract, although its value may differ at
each occurrence.
In fact, the second line of $(2.6)$ implies readily that for $\sigma>1$
$$
|L(s;\r{sym}^\ell V\times \r{sym}^\ell V)|^{-1}
\le L(\sigma;\r{sym}^\ell V\times \r{sym}^\ell V).\eqno(2.15)
$$
To the right side we apply $(5.6)$ below 
if $\ell=2$ and an analogous bound
if $\ell=1$ that is obtainable in much the same way as $(5.6)$.
This yields $(2.14)$. We then apply either [13, (1.4.31)] or
[19, Lemma $\alpha$] to the function $f(s)=(s-1)
L(s;\r{sym}^\ell V\times \r{sym}^\ell V)$ 
which is entire and of polynomial growth in both $s$ and $\nu_V$ 
whenever $\Re s$ is bounded, as is to be explained below.
We get, for $s_0=1+a\big(\log(|\nu_V|+|\Im s|)\big)^{-1}$,
$$
\bigg|{f'\over f}(s)-\sum_\rho{1\over s-\rho}\bigg|\ll 
\log (|\nu_V|+|\Im s|),\quad |s-s_0|\le C,
\eqno(2.16)
$$
where $C>0$ is any constant larger than 10, say, and
$\rho$ runs over zeros of $f(s)$ 
such that $|\rho-s_0|\le C/2$. Specialising $(2.16)$
with $s=\sigma_0=1+a(\log|\nu_V|)^{-1}$, it follows that
$$
-{L'\over L}(\sigma_0;\r{sym}^\ell V\times\r{sym}^\ell V)
={1\over \sigma_0-1}-\sum_\rho{1\over\sigma_0-\rho}
+O(\log|\nu_V|),\eqno(2.17)
$$
where the sum over $\rho$ is real and non-negative, since
$f(\overline{\rho})=0$. On the
other hand the middle line of $(2.6)$ implies that the left side is
also non-negative. Thus, we have, for $\ell=1,2$,
$$
{L'\over L}(\sigma_0;\r{sym}^\ell V\times\r{sym}^\ell V)
\ll(\sigma_0-1)^{-1}+\log|\nu_V|.\eqno(2.18)
$$
The same bound holds for $(f'/f)(\sigma_0)$ obviously.
Then we appeal to either [13, (1.4.34)] or [19, Lemma $\gamma$],
and find in particular that $(f'/f)(1)\ll\log|\nu_V|$. Combining this bound
and $(2.10)$, we obtain $(2.13)$.
\par
It would be remiss not to make it explicit that the $L$-functions
dealt with in the above, save for $(2.10)$ with
general $\ell$, are known to be of finite order
with respect to $s$. This follows either from the
cuspidality of the relevant representations or from
the integral representations for Rankin convolutions
via Eisenstein series.  It endorses the application of 
the Phragm\'en--Lindel\"of convexity principle and
leads us to the polynomial growth of the $L$-functions which
is uniform even in $\nu_V$.
Thus, for our immediate purpose, we are assured that
for $|\Re s|\le 2$ and $\ell=1,2$,
$$
\eqalignno{
&(s-1)L\big(s;\r{sym}^\ell V\times\r{sym}^\ell V\big)
\ll (|\nu_V|+|s|)^c, &(2.19)\cr
&\quad\qquad L\big(s;\r{sym}^{2\ell}V\big)
\ll (|\nu_V|+|s|)^c,&(2.20)
}
$$
upon which essentially the
whole of our analysis is laid, including basic results
such as $(2.12)$ as well. We should add also that
the same finiteness assertion for $L$-functions
within a broader framework is 
established in Gelbart--Shahidi [3]; it appears to us
that prior to their work the fact had been assumed oft-tacitly
or rather treated as an assertion too obvious to mention explicitly.
This ends the proof of Lemma 1.
\medskip
\noindent
{\bf Lemma 2.} {\it The number of zeros of $L(s;\r{sym}^{2\ell}V)$,
$\ell=1,2$, on the disk of radius $\tau$ with centre 
at $1+it,\, t\in\B{R}$, is 
$O\big(\tau\log(|\tau_V|+|t|)\big)$, providing
$a(\log(|\nu_V|+|t|))^{-1}\le\tau\le C/8$, with
$a$, $C$ as in $(2.12)$ and $(2.16)$, respectively.\/}
\par
\noindent
{\it Proof\/}.  It suffices to prove the same assertion on the function
$f$ utilised in the above. We apply $(2.16)$ with $s=1+\tau+it$. We get
$$
\Re\sum_\rho{1\over s-\rho}\ll\log (|\nu_V|+|t|),\eqno(2.21)
$$
since $|(f'/f)(s)|\le|(f'/f)(1+\tau)|\ll \log(|\nu_V|+|t|)$ by 
just the same reasoning leading $(2.18)$.
Restricting the sum to those zeros such that
$|s-\rho|\le 2\tau$, we end the proof.
\bigskip
\noindent
{\bf 3. Sieve tools}
\par
\noindent
Now, the argument in [11][12][13, Vol.\ I, Section 9.3]
yields readily that the optimal
choice of $\{\lambda_d\}$ for $(2.9)$ is given by
$$
\lambda_d=\mu(d)F_d{G_d(R/d)\over G_1(R)};
\eqno(3.1)
$$
in fact, this could rather be set a priori without any sieve reasoning. 
Here $\mu$ is the M\"obius function,
$$
F_d=\prod_{p|d}F_p,\quad
G_d(x)=\sum_{\scr{r\le x}\atop\scr{(d,r)=1}}
{\mu^2(r)K(r)},\quad K(r)=\prod_{p|r}
(F_p-1),\eqno(3.2)
$$
where $(d,r)$ is
the greatest common divisor of $d$ and $r$,
and $F_p=F_p(1)$ with
$$
F_p(s)=\sum_{l=0}^\infty \tau^{(2)}_{V\times V}
(p^l)p^{-l s}\eqno(3.3)
$$
which is the $p$-factor of $(2.6)$, $\ell=2$, 
converging absolutely for $\Re s>{7\over16}$
because of $(1.5)$.
It should be stressed that for any $p$
$$
F_p-1\ge {1\over p^3}.\eqno(3.4)
$$
In fact, we have, by the second line of $(2.6)$,
$$
\eqalignno{
F_p-1&\,\ge\big(\alpha_V(p)^2+1
+\alpha_V(p)^{-2}\big)^2{1\over p}
+\big(\alpha_V(p)^6+1+\alpha_V(p)^{-6}\big)^2
{1\over3p^3}\cr
&\,=\tau_V^2(p^2){1\over p}
+\big(\tau^3_V(p^2)-
3\tau^2_V(p^2)+3\big)^2{1\over3p^3}.&(3.5)
}
$$
If $|\tau_V(p^2)|\ge{1\over2}$, then $(3.4)$ is obvious;
otherwise the multiplier of $1/3p^3$ is larger than 
$4$.

\medskip
The choice $(3.1)$ 
leads us to the multiplicative function $\Phi_r$:
$$
\eqalign{
&\,\Phi_r(n)={\mu((r,n))\over K((r,n))},\cr
\sum_{d|n}\lambda_d=&\,{1\over G_1(R)}\sum_{r\le R}
\mu^2(r)K(r)\Phi_r(n);
}\eqno(3.6)
$$
see [11][12, \S1.4][13, Vol.\ I, Chapter 9] as well as (iii) of
Section 5. We are about to show the
quasi-orthogonality in the set $\{\Phi_r(n): n\in\B{N}\}$.
To this end, we consider the expression
$$
\sum_{N\le n\le M+N}\tau^{(2)}_{V\times V}(n)
\Bigg|\sum_{r\le R}
\mu^2(r)\Phi_r(n)\sqrt{K(r)}\cdot
b_r\Bigg|^2,\eqno(3.7)
$$
where $M,N, R\ge1$, and $\{b_r\}$ are all arbitrary.
Expanding the squares out, we have
$$
\sum_{r_1,\,r_2\le R}{\mu^2(r_1)\mu^2(r_2)
\sqrt{K(r_1)K(r_2)}}\big\{S(M+N;r_1,r_2)
-S(N;r_1,r_2)\big\}
b_{r_1}\overline{b}_{r_2},\eqno(3.8)
$$
with
$$
S(N; r_1,r_2)=
\sum_{n\le N}\tau_{V\times V}^{(2)}(n)
\Phi_{r_1}(n)\Phi_{r_2}(n).\eqno(3.9)
$$
Thus, let us consider the function
$$
\eqalignno{
&\sum_{n=1}^\infty \tau^{(2)}_{V\times V}(n)
\Phi_{r_1}(n)\Phi_{r_2}(n)n^{-s}\cr
=&\,\Bigg(\sum_{(n,\, r_1r_2)=1}\Bigg)
\Bigg(\sum_{n|([r_1,\,r_2]/(r_1,\,r_2))^\infty}\Bigg)
\Bigg(\sum_{n|( r_1,\,r_2)^\infty}\Bigg)
=\r{F}_1\r{F}_2\r{F}_3,&(3.10)
}
$$
say, where $[r_1,r_2]$ is the least common multiple of
$r_1$ and $r_2$, and it is assumed temporarily that $\Re s$ is
sufficiently large. We have
$$
\eqalign{
\r{F}_1&=L(s;\r{sym}^2V\times\r{sym}^2V)
\prod_{p\mid r_1r_2}F_p(s)^{-1},\cr
\r{F}_2&=\prod_{p|{\scr{[r_1,\,r_2]}\over\scr{(r_1,\,r_2)}}}
\big(1-(F_p-1)^{-1}(F_p(s)-1)\big),\cr
\r{F}_3&=\prod_{p|(r_1,\,r_2)}
\big(1+(F_p-1)^{-2}(F_p(s)-1)\big).
}\eqno(3.11)
$$
We write this as
$$
\r{F}_1\r{F}_2\r{F}_3=L(s;\r{sym}^2V\times
\r{sym}^2V)U_{r_1,r_2}(s);\eqno(3.12)
$$
thus
$$
S(N;r_1,r_2)=
\sum_{d|(r_1r_2)^\infty} u(d)
\sum_{n\le N/d}\tau_{V\times V}^{(2)}(n),\eqno(3.13)
$$
where $u(d)$ is the coefficient of the Dirichlet series
$U_{r_1,r_2}(s)$ and empty sums are to vanish. 
\par
To the last inner sum we apply the asymptotic formula
$$
\sum_{n\le N}\tau^{(2)}_{V\times V}(n)
=L\big(1;\r{sym}^2V\big)L\big(1;\r{sym}^4V\big)N
+O\big(|\nu_V|^c N^{4/5+\eta}\big),
\eqno(3.14)
$$
with any fixed $\eta>0$; see (v) of
Section 5. We now have that
$$
S(N;r_1,r_2)={\delta_{r_1,r_2}\over K(r_1)}
L(1;\r{sym}^2V)L(1;\r{sym}^4V)N
+O\big(|\nu_V|^c (r_1r_2)^{7/2}N^{5/6}\big),
\eqno(3.15)
$$
with the Kronecker delta. This is due to the facts that
$U_{r_1,r_2}(1)=\delta_{r_1,r_2}/K(r_1)$ by $(3.11)$
and that a combination of $(3.4)$, $(3.11)$ and
the expression for $F_p(s)^{-1}$ inferred from the
first line of $(2.6)$ gives
$$
\sum_{d|(r_1r_2)^\infty}|u(d)|d^{-5/6}\le
c^{\nu(r_1r_2)}(r_1r_2)^3,\eqno(3.16)
$$
with $\nu(m)=\sum_{p|m}1$. In fact, the left side is
not greater than
$$
\eqalignno{
&\,\prod_{p|r_1r_2}\prod_{j=0}^2\prod_{k=0}^2
\Big(1+|\alpha_V^{2(2-j-k)}(p)|p^{-5/6}\Big)\cr
\times&\,\prod_{p|{\scr{[r_1,\,r_2]}\over\scr{(r_1,\,r_2)}}}
\Big(1+(F_p-1)^{-1}\big(F_p\big(\txt{5\over6}\big)
-1\big)\Big),\cr
\times&\,\prod_{p|(r_1,\,r_2)}
\Big(1+(F_p-1)^{-2}\big(F_p\big(\txt{5\over6}\big)
-1\big)\Big)\cr
\le&\; c^{\nu(r_1r_2)}
\prod_{p|r_1r_2}(F_p-1)^{-1},&(3.17)
}
$$
in which we have applied $(1.5)$.
One may show a better bound, but $(3.16)$ suffices
for our purpose. 
\medskip
Collecting these assertions and
invoking the duality principle together with $(2.5)$, 
we obtain the
following analogue of [12, Theorem 5]
[13, Vol.\ I, (9.1.23)]:
\medskip
\noindent
{\bf Lemma 3.} {\it We have, uniformly
for $1\le M\le N$, $R\ge1$ and for any complex vector
$\{a_n\}$, 
$$
\eqalignno{
&\sum_{r\le R}\mu^2(r)K(r)
\Bigg|\sum_{N\le n\le M+N}\tau_V^{(2)}(n)
\Phi_r(n)a_n\Bigg|^2\cr
\ll &\,\Big(L(1;\r{sym}^2V)L(1;\r{sym}^4V)M
+|\nu_V|^cR^8N^{5/6}\Big)
\sum_{N\le n\le M+N}|a_n|^2.&(3.18)
}
$$
\/}
\medskip
\noindent
As a corollary, we have
\medskip
\noindent
{\bf Lemma 4.} {\it Let ${\cal S}=\{s_j\}$ 
be a finite set of complex numbers such that 
$\Re s_j\ge0$, $|\Im s_j|\le T$ and $|\Im(s_j-s_k)|
\ge\xi>0$, $j\ne k$. Then we have, for any complex
sequence $\{a_n\}$,
$$
\eqalignno{
&\sum_{r\le R}\mu^2(r)K(r)\sum_{s\in{\cal S}}
\Bigg|\sum_{n=1}^\infty\tau_V^{(2)}(n)
\Phi_r(n)a_nn^{-s}\Bigg|^2\cr
\ll &\,\sum_{n=1}^\infty(\xi^{-1}+\log n)
\Big(L(1;\r{sym}^2V)L(1;\r{sym}^4V)n
+|\nu_V|^cR^8Tn^{5/6}\Big)
|a_n|^2,&(3.19)
}
$$
provided the right side converges.\/}
\smallskip
\noindent
{\it Proof\/}. This is essentially a specialisation of [13, 
Vo.\ I, (9.1.42)] which
originates in [10, (7.7)]
and [12, Lemma 26] (the corresponding
display line in its digitised edition published by Tata IFR 
contains an obvious misprint).
\medskip
On the other hand, the sieve effect of $(3.1)$ is embodied in 
\medskip
\noindent
{\bf Lemma 5.} {\it Provided 
$$
\hbox{$\log R/\log|\nu_V|$ is sufficiently large but bounded,}
\eqno(3.20)
$$
we have
$$
G_1(R)\asymp L(1;\r{sym}^2V)
L(1;\r{sym}^4V)\log R.\eqno(3.21)
$$
In particular, we have
$$
L(1;\r{sym}^2V)L(1;\r{sym}^4V)\gg 
(\log|\nu_V|)^{-1}.\eqno(3.22)
$$
}
\par
\noindent
{\it Proof\/}. We have, in the region of absolute
convergence,
$$
\sum_{r=1}^\infty\mu^2(r)K(r)r^{-s}=
L(s+1;\r{sym}^2V\times\r{sym}^2V)Y(s),\eqno(3.23)
$$
where
$$
Y(s)=\prod_p\Bigg\{\Big(1+(F_p-1)p^{-s}\Big)
\prod_{j=0}^2\prod_{k=0}^2
\Big(1-\alpha_V^{2(2-j-k)}(p)p^{-s-1}\Big)\Bigg\}
\eqno(3.24)
$$
which is absolutely convergent and bounded for $\Re s>-
{1\over 16}$ because of $(1.5)$. 
We multiply both sides of $(3.23)$
by $\Gamma(s)R^s/2\pi i$ and 
integrate over the line $\Re s=1$. The
left side of the resulting identity
is obviously $\gg G_1(R)$. On the right side, we
may shift the contour to
$\Re s=-{1\over 17}$, encountering only one singularity
which is a double pole at the origin, because of
$(2.10)$, $\ell=2$, and Lemma 1; note that $(2.19)$ is
necessary here as well. The residue equals 
$$
\eqalignno{
&L(1;\r{sym}^2V)
L(1;\r{sym}^4V)\log R\cr
\times&\,\bigg[
1+O\Big\{\Big(\big|(L'/L)(1;\r{sym}^2V)\big|
+\big|(L'/L)(1;\r{sym}^4V)\big|\Big)/\log R\Big\}
\bigg], &(3.25)
}
$$
as $Y(0)=1$ and $Y'(0)$ is bounded. The new integral
is negligible, since we have imposed $(3.20)$ and the
left side of the integrated identity 
is greater than or equal to $\exp(-1/R)$ 
by definition.
As to the logarithmic derivatives in $(3.25)$, we apply
$(2.13)$. This proves the upper bound for $G_1(R)$
contained in $(3.21)$. To prove lower bound, we
multiply $(3.23)$ by $\Gamma(s)(R/\log^2\!R)^s/2\pi i$
and proceed in just the same way. Lastly, 
to prove $(3.22)$ it should
suffice to observe that $G_1(R)>1$. This ends the proof.
\medskip
Further, we need to invoke [12, Theorem 4][13, Vol.\ I,
(9.2.11)]:
\medskip
\noindent
{\bf Lemma 6.} {\it Let $v$ be a large positive parameter,
and let $\vartheta>0$ be a constant. 
We put, with an integer $l\ge0$,
$$
\displaystyle{\Xi_d^{(l)}={1\over l!}(\vartheta\log v)^{-l}
\sum_{j=0}^l(-1)^{l-j}{l\choose j}\xi_d^{(j,l)},}\atop
\hbox{$\xi_d^{(j,l)}=\mu(d)
\Big({\log v^{1+j\vartheta}/d\Big)^l} $ 
for $d\le v^{1+j\vartheta}$, and  $= 0$ 
 for $d> v^{1+j\vartheta}$}.\eqno(3.26)
$$
Then, we have that
$$
\hbox{$\Xi_d^{(l)} =\mu(d)$ for $d\le v$}\eqno(3.27)
$$
and 
$$
\sum_{n=1}^\infty d_l(n)
\bigg(\sum_{d|n}\Xi_d^{(l)}\bigg)^2n^{-\omega}\ll 1,
\eqno(3.28)
$$
whenever $ \omega\ge 1+1/\log v$.
Here $d_l(n)$ is the number of ways of representing 
$n$ as a product of $l$ positive integral factors, and the
implied constant depends on $l$ and $\vartheta$ 
at most.\/}
\bigskip
\noindent
{\bf 4. Proofs of theorems}
\smallskip
\noindent
{\it Proof of Theorem 2\/}.
With Lemmas 1--6 in hand, the discussion is essentially 
a repetition of the argument
developed in [12, \S5.2]. Thus, in order to estimate
$N_{\r{sym}^2V}(\alpha,T)$, we may assume obviously
that ${4\over5}\le\alpha\le1-{a/\log(|\nu_V|T)}$,
where $a$ is as in $(1.11)$.
For each $\kappa=0,1$ we pick up
a zero of $L(s;\r{sym}^2V)$ lying simultaneously
in the rectangle in question and in one of the horizontal
strips $(2n+\kappa)/\log T\le \Im s 
<(2n+\kappa+1)/\log T$, $n\in\B{Z}$;
the resulting set of zeros is denoted by $\e{Z}_\kappa$.
Since Lemma 2 implies that the number of zeros in the disk $|s-(1+iu)|\le 
1-\alpha$ is $\ll (1-\alpha)\log(|\nu_V| T)$ 
whenever $-T\le u\le T$, we have
$$
N_{\r{sym}^2V}(\alpha,T)\ll (1-\alpha)(\log |\nu_V|T)
\Big\{|\e{Z}_0|+|\e{Z}_1|\Big\}.\eqno(4.1)
$$
Then, we make a conversion of [12, Lemma 5]
[13, Vol.\ I, (9.2.2)]
to our present situation: with
$\mu^2(r)=1$, we have, in view of $(3.27)$,
$$
\eqalignno{
&\qquad 1+\sum_{v\le n} \tau_V^{(2)}(n)\Phi_r(n)
\bigg(\sum_{d|n}\Xi^{(1)}_d\bigg)n^{-s}
=L(s;\r{sym}^2V)J_r(s),&(4.2)\cr
&J_r(s)={1\over K(r)}\sum_{d=1}^\infty 
\mu((d,r))\Xi^{(1)}_d
\prod_{p|d}\big(1-X_p(s)\big)
\prod_{\scr{p\nmid d}\atop\scr{p|r}}
\big(X_p(s)F_p-1\big),&(4.3)
}
$$
for $\Re s>1$, where $F_p$ is 
as in the previous section and
$X_p(s)$ is the inverse of the $p$-factor of
$(2.2)$, $\ell=2$. We set in Lemmas 3--6
$$
R=(|\nu_V|T)^A,\; v=R^A,\; \vartheta=1/A,\eqno(4.4)
$$
with a sufficiently large constant $A>0$. 
Then, $J_r(s)$ is entire, and $J_r(s)\ll v$,
for $r\le R$ and $\Re s\ge{3\over4}$, which can be confirmed
readily by applying $(1.5)$ and $(3.4)$.
With this, let us consider the expression
$$
\int_{2-i\infty}^{2+i\infty}
L(\rho+w; \r{sym}^2V)J_r(\rho+w)
\Gamma(w)v^{2w} dw,\quad \rho\in\e{Z}_\kappa.
\eqno(4.5)
$$
Invoking $(2.15)$, $\ell=2$, again and 
shifting the contour to the left appropriately, we see that
we do not encounter any singularity and
$(4.5)$ is negligibly small; hence,
$$
{1\over2}\le\Bigg|\sum_{v\le n\le  
v^4}\tau^{(2)}_V(n)\Phi_r(n)
\bigg(\sum_{d|n}\Xi_d^{(1)}\bigg)n^{-\rho}
e^{-n/v^2}\Bigg|^2.\eqno(4.6)
$$
We multiply both sides by the factor 
$\mu^2(r)K(r)$ and sum over $r\le R$ as well as over
$\rho\in\e{Z}_\kappa$, getting
$$
\eqalignno{
&G_1(R)\big|\e{Z}_\kappa\big|\ll 
\sum_{r\le R}\mu^2(r)K(r)\cr
&\times\sum_{\rho\in\e{Z}_\kappa}
\Bigg|\sum_{v\le n\le  
v^4}\tau^{(2)}_V(n)\Phi_r(n)
\bigg(\sum_{d|n}\Xi_d^{(1)}\bigg)n^{-\rho}
e^{-n/v^2}\Bigg|^2.&(4.7)
}
$$
By virtue of Lemma 4 and $(3.22)$, we have that
$$
\eqalignno{
G_1(R)\big|\e{Z}_\kappa\big|
\ll&\, L(1;\r{sym}^2V)L(1;\r{sym}^4V)(\log|\nu_V|T)
v^{4(1-\alpha)}\cr
&\times\,\sum_{n=1}^\infty \bigg(\sum_{d|n}\Xi_d^{(1)}
\bigg)^2n^{-\omega_0},&(4.8)
}
$$
with $\omega_0=1+(\log |\nu_V|T)^{-1}$.
Then,  by $(3.28)$ with $l=1$, we find that
$$
G_1(R)\big|\e{Z}_\kappa\big|\ll
 L(1;\r{sym}^2V)L(1;\r{sym}^4V)(\log|\nu_V|T)
v^{4(1-\alpha)}.\eqno(4.9)
$$
Therefore, in view of $(3.21)$, 
we end the proof of Theorem 2.
\medskip
\noindent
{\it Proof of Theorem 1.\/} Let $N$ be a parameter tending to
infinity. Let $g$ be a $C^\infty$-function 
such that $g(x)=1$ for $N-M\le x\le N$; $0\le g(x)\le1$ for
either $N-M-U\le x\le N-M$ or $N\le x\le N+U$;
$g(x)=0$, otherwise. Here 
$M=N^{1-\theta}$ and $U\le M/2$,
with a small parameter 
$\theta>0$. We have, provided $U\ge N^{2/3}$,
$$
\sum_{N-M\le p\le N}(\log p)\tau_V^2(p)
=\sum_p (\log p)\tau_V^2(p)g(p)
+O\big(|\nu_V|^b U\log N\big),\eqno(4.10)
$$
with a constant $b>0$. This lower bound for $U$
comes from the asymptotic formula for  
$\sum_{n\le N}\tau_V^2(n)$
which can be inferred from the discussion 
in (v) of the next section:
it corresponds to the case $r=4$ there. 
Then, we set $U=N^{1-(b+2)\theta}\ge N^{2/3}$, 
which induces an upper bound for
$\theta$. We have $|\nu_V|^bU\le MN^{-\theta}$, 
provided $|\nu_V|\le N^\theta$. 
\par
Next, let $\hat{g}$
be the Mellin transform of $g$; then we have, in $(4.10)$,
$$
g(p)={1\over 2\pi i}\int_{\xi-iQ}^{\xi+iQ}\hat{g}(s)p^{-s}ds
+ O\big(N^{-A\theta}(N/p)^\xi\big),\eqno(4.11)
$$
where $Q=N^{1+\theta}/U$,
with arbitrary $\xi,\,A>0$ on which the implied constant may
depend.
This is due to the fact that 
a multiple application of integration by parts
gives $\hat{g}(s)\ll (N/U)^{k-1}N^{\Re s}|s|^{-k}$ 
for each $k\in\B{N}$,
as one may assume that $g^{(k)}(x)\ll U^{-k}$;
see [13, Vol.\ II, p.\ 41], for instance. Before
inserting $(4.11)$ into $(4.10)$, we observe that 
by $(2.1)$
$$
\eqalignno{
-{L'\over L}\big(s;V\times V\big)&\,
=-{\zeta'\over\zeta}(s)-
{L'\over L}\big(s;\r{sym}^2V\big)\cr
&\,=\sum_p(\log p)\tau_V^2(p)p^{-s}+E(s;V),&(4.12)
}
$$
for $\Re s>1$, where $E(s;V)$ is regular and bounded for 
$\Re s\ge{3\over4}$ because of $(1.5)$.
Hence, we now have, with $\xi=1+(\log N)^{-1}$,
$$
\sum_{N-M\le p\le N}(\log p)\tau_V^2(p)=
-{1\over2\pi i}\int_{\xi-iQ}^{\xi+iQ}{L'\over L}(s; V\times V)
\hat{g}(s)ds +O\big(MN^{-\theta}\log N\big).\eqno(4.13)
$$
In fact, the
contribution of $E(s;V)$ can be seen to be negligible
after moving the contour to $\Re s={3\over4}$ and noting 
$\hat{g}(s)\ll N^{3/4}|s|^{-1}$ there; that of the
error term in $(4.11)$ is $O\big(N^{1-A\theta}(1+
|(L'/L)(\xi; V\times V)|)\big)$ and this logarithmic derivative
is $O\big(\log(|\nu_V|N)\big)$, similarly to
$(2.18)$.
\par
We may adjust $Q$ so that $(L'/L)(\sigma+iQ;V\times V)\ll 
(\log |\nu_V|Q)^2$ for $-1\le\sigma\le 2$, which can
be shown by using $(2.16)$ and Lemma 2 via a well-known
argument (an application of the pigeon box principle).
We combine this fact 
with the above estimate of $\hat{g}(s)$ and find that
$$
\sum_{N-M\le p\le N}(\log p)\tau_V^2(p)=M-
\sum_\rho\hat{g}(\rho)+O\big(MN^{-\theta}\log N\big),
\eqno(4.14)
$$
where $\rho$ runs over all the 
zeros of $L(s;V\times V)$ in the
part with $|\Im s|\le Q$ of the critical strip; the details of this
procedure are skipped, since they are analogous to dealing with
$\zeta(s)$ instead.
We have, by $(1.11)$--$(1.12)$,
$$
\eqalignno{
\sum_\rho|\hat{g}(\rho)|&\,\le\int_{N-M-U}^{N+U}
g(x)(\log x)\Bigg(\int_0^1 N_{V\times V}(\alpha,Q)
x^{\alpha-1}d\alpha\Bigg) dx\cr
&\,\ll M(\log N)\int_0^{1-a/\log(|\nu_V|Q|)}
\Big((|\nu_V|Q)^\omega N^{-1}\Big)^{1-\alpha}
d\alpha.&(4.15)
}
$$
Here we set $\theta\le 1/(2(b+4)\omega)$ so that
$(|\nu_V|Q)^\omega\le N^{1/2}$; then, 
the choice of the lower limit of the integral over $\alpha$ is
irrelevant. Hence, under this
restriction on $\theta$, we find that
$$
\sum_\varrho|\hat{g}(\varrho)|\ll M \exp\Big(
-{a\over 2(b+4)\theta}\Big),\eqno(4.16)
$$
where the implied constant is absolute, especially
independent of $\theta$. Finally, we
note that if $(\log N)^{-1/2}\le \theta$, then
$N^{-\theta}\log N<\exp(-1/2\theta)$; 
namely, if $a$ is adjusted to be small enough, then
the error term in $(4.14)$ is
absorbed in the right side of $(4.16)$. We still need to
eliminate the factor $\log p$ on the left side $(4.14)$.
It should, however, suffice to observe that 
$\log p=\log N +O(M/N)$, and this error term contributes
$O\big(M^2/(N\log N)\big)$, since $(4.14)$ itself is $O(M)$,
as it follows from what we have discussed so far.
We end the proof of Theorem 1.
\bigskip
\noindent
{\bf 5. Concluding Remarks} 
\smallskip
\noindent
(i) Hoheisel's prime number theorem [5] asserts that there exists 
a constant $\theta_1>0$ such that
$$
\sum_{x-y\le p\le x}1=\big(1+o(1)\big)
{y\over \log x},\quad y=x^{1-\theta_1},
\eqno(5.1)
$$
as $x$ tends to infinity. Thus our theorem is in fact
a {\it partial\/} analogue of Hoheisel's theorem in the sense that 
both $(1.10)$ and $(5.1)$ imply the existence 
of primes in short intervals $[x-y, x]$, 
while $(1.10)$ is not an asymptotic
identity if $\theta$ there remains independent of $x$.
Essential ingredients of Hoheisel's argument are
the zero-free region
 $$
\zeta(s)\ne0,\quad \Re s>
1-c{\log\log(|\Im s|+3)\over\log(|\Im s|+3)},\eqno(5.2)
$$
due to Littlewood, and the zero-density estimate
$$
N(\alpha,T)\ll T^{4\alpha(1-\alpha)}(\log T)^6,\eqno(5.3)
$$
due to Hoheisel himself,
with $N(\alpha,T)$ being the number of zeros of $\zeta(s)$
in the rectangle ${1\over2}\le\alpha\le\Re s\le1$, 
$|\Im s|\le T$,
$T\ge2$. Hoheisel argued in a manner quite
similar to how we have argued on $(4.15)$, and he
needed $(5.2)$ in order to offset the presence of the
logarithmic factor in $(5.3)$. Later a 
zero-free region far superior
than $(5.2)$ was established by I.M. Vinogradov, and 
various strengthening of $(5.3)$
followed, resulting in 
numerical improvements upon the constant $\theta_1$
in $(5.1)$. Nevertheless, the overall structure of
the zero-density method in the theory of
the distribution of primes has remained essentially the same
to this day
since Hoheisel's pioneering work; see [10] [13, Vol.\ I] for
the relevant history. In this context,  
our extension of $(5.1)$ to Maass forms, even though it is
only partial, is a more delicate work, since
any analogue of $(5.2)$ or the like is unknown
in the entire theory of $L$-functions associated with cusp forms; 
indeed, this remains a challenging open problem. 
In order to compensate for
this difficulty, one needs the zero-density estimate $(1.12)$. 
The entire argument in the present article
is devoted to the elimination of the logarithmic
factor that could come up in bounding 
$N_{V\times V}(\alpha,T)$. 
\par
We are, incidentally,
extremely grateful to A. Perelli for sending
us a photo-copy of [5] which would have been hard for us
to access otherwise. Knowing of our quest for a copy
of the paper, he kindly made a search through 
his impressive collection of off-prints. Included in the collection
are some complimentary copies, originally presented
to R. Rankin; it was amongst these Hoheisel's paper was found. 
Rankin had given M. Nair boxes of off-prints, 
from which Nair made
a gift set for his friend. This episode is a pleasant
coincidence, for it further signifies that our work rests upon
a combination of Hoheisel's and Rankin's fundamental ideas. 
\medskip
\noindent
(ii) It might be worth remarking that one may avoid 
appealing to $(1.11)$, since it is possible to prove an analogue
of the Deuring--Heilbronn--Linnik phenomenon for the function
$L(s;V\times V)$ itself, by following
the argument of [15]. Namely, 
one may allow the possibility
of the existence of an exceptional zero 
for $L(s;V\times V)$, as
it should repel all other zeros toward the left of
the line $\Re s=1$ deeper than $(1.11)$.
\medskip
\noindent
(iii) As to $(3.6)$, the use of quasi-characters in 
the study of the zero density of Dirichlet $L$-functions was
initiated by Selberg [18]; actually the quasi-orthogonality
among Ramanujan sums was indicated there. 
However, his character does not straightforwardly
generalise to $(3.6)$. To identify ours,
we need the observation made in [11, p.\ 166] [12, p.\ 40]
 [13, Vol.\ I, Chapter 9] 
that Selberg's quasi-character originates
in fact in the $\Lambda^2$-sieve applied to the most primitive
arithmetic function: the constant $1$. With
this, one may come to the idea $(2.4)$, and to 
the definition $(2.9)$. Nonetheless, the employment of $(2.9)$ 
in place of $(2.4)$ is never trivial itself. The use of the factor 
$\tau^{(2)}_{V\times V}(n)$ is
deliberately made in order to exploit the Rankin convolution, 
which is
not necessary when dealing with holomorphic cusp forms
because of the validity of the Ramanujan conjecture; 
see [15, Part I]. We
stress that our argument generalises to the treatment of
sums of Hecke--Maass eigenvalues raised to the 
$2\ell$-th power
over primes in short intervals, provided the holomorphy, 
the polynomial 
growth and the zero-free region of the type $(2.12)$ are all
available for the quotient
$L(s;\r{sym}^{2\ell}V\times\r{sym}^{2\ell}V)/\zeta(s)$.
\medskip
\noindent
(iv) In [16] we shall discuss an analogue of
Linnik's least prime number theorem for Maass forms; namely,
an extension of Theorem 1 to arithmetic
progressions is to be achieved.
There we shall need to develop a large sieve zero density
estimate as well as an analogue of the
Linnik phenomenon for the $\chi$-twist of 
$L(s;\r{sym}^2V)$, which can be treated
in a combined fashion as an extension of Lemma 3. 
This is an interesting issue in itself.
\medskip
\noindent
(v) The asymptotic formula $(3.14)$ is an essential 
ingredient in our argument. It
might appear to be a corollary of a basic result 
due to Chandrasekharan--Narasimhan 
[1, Theorem 4.1] which itself stems from Landau [8].
However, their assertion is not uniform, 
especially with respect to the shift parameters involved in
their $\Gamma$-factors: $\nu_V$ in our case.
Hence, we take here the task to give 
a brief proof of $(3.14)$. 
Despite the present specialisation,
our argument should extend to the general situation
that is discussed in [1, Section 4]. Ours is, in fact,
based on the treatment of the
divisor problem developed in [13, Vol.\ II, Section 1.1].
Note that $(2.19)$ is indispensable hereafter.
\par
Thus, let $w$ be a test function compactly
supported on the positive real axis. Let $\hat w$ be its Mellin
transform. We have, by a standard procedure,
$$
\eqalignno{
\sum_{n=1}^\infty\tau_{V\times V}^{(2)}(n)w(n)
&\,=L(1;\r{sym}^2V)L(1;\r{sym}^4V)\hat{w}(1)\cr
&\,+\kappa_{V,r}\sum_{n=1}^\infty
\tau_{V\times V}^{(2)}(n)W_V(n)/n,&(5.4)
}
$$
with $|\kappa_{V,r}|=1$ and
$$
W_V(n)={1\over2\pi i}
\int_{-{1\over 2}-i\infty}^{-{1\over2}+i\infty}
{\overline{\Gamma_r(1-\bar{s})}
\over\Gamma_r(s)}\hat{w}(s)n^sds.\eqno(5.5)
$$
Here it is understood that the $\Gamma$-factor 
of the functional equation for $L(s;\r{sym}^2V\times
\r{sym}^2V)$ is denoted as $\Gamma_r(s)$; thus
$r=9$ in fact; and one may assume without any loss of
generality that $\gamma_V^{(r,j)}$ in $(1.7)$
are such that we do not encounter any singularity
save for the simple pole at $s=1$ while
shifting the contour. 
We observe that with an appropriate choice of $w$ the
expansion $(5.4)$ implies in particular
$$
\sum_{n\le N}\tau_{V\times V}^{(2)}(n)\ll |\nu_V|^cN,
\eqno(5.6)
$$
uniformly in $N\ge1$.
This is due to the fact that the residue of $L(s;\r{sym}^2V
\times\r{sym}^2V)$ at $s=1$ as well as the size of
the function itself on $\Re s={3\over2}$ are
$\ll|\nu_V|^c$, as can be seen  via a combination of
the non-negativity of $\tau^{(2)}_{V\times V}(n)$,
the functional equation and $(2.14)$.
On the other hand, we have
$$
W_V(n)=\int_0^\infty w^{(q)}(x)x^{q-1}
G_V(nx;q)dx,\eqno(5.7)
$$
where 
$$
G_V(y;q)={(-1)^q\over2\pi}\int_{-\infty}^\infty
{\overline{\Gamma_r(1-\bar{s})}\over s(s+1)
\cdots(s+q-1)\Gamma_r(s)}y^sdt,\quad
s=-\txt{1\over2}+it,\eqno(5.8)
$$
with any sufficiently large $q$.
We assume that 
$$
\hbox{$\log y/\log |\nu_V|$ is large.}\eqno(5.9)
$$ 
By the Stirling formula,
we see that the saddle point of the last integral is
at $t=t_0=2\pi y^{1/r}$. 
With this in mind, we set $G_V=G_V^{(1)}+G_V^{(2)}$, 
in which the first term on the right
corresponds to $|t-t_0|<{1\over2}t_0$ and
the second to the rest; one may suppose that
smooth weights $\eta_1(t)$ and $\eta_2(t)$ 
have been attached 
to the integrands of $G_V^{(1)}$ and 
$G_V^{(2)}$, respectively, 
in an obvious manner. Accordingly,
we have $W_V=W_V^{(1)}+W_V^{(2)}$.
Then the off-saddle 
contribution $W_V^{(2)}$ turns out
to be negligible via a multiple application 
of integration by parts with respect to $t$ in $G_V^{(2)}$. 
As to $W_V^{(1)}$,
we undo the integration by parts with respect
to $x$ in $(5.7)$, going back to $q=0$.
That is, we have
$$
W_V^{(1)}(n)=\int_0^\infty w(x)x^{-1}G_V^{(1)}
(nx)dx,\eqno(5.10)
$$
with $G_V^{(1)}(y)=G_V^{(1)}(y;0)$.
We get, after some computation by means of Stirling's formula,
$$
G_V^{(1)}(y)={1\over2\pi\sqrt{y}}
\int_{-\infty}^\infty \eta_1(t)\left[{t\over2\pi}
\exp\left(-it\log{t\over2\pi ey^{1/r}}\right)\right]^r
B_1(t)dt,\eqno(5.11)
$$
where $B_1(t)$ is a power series in $1/t$, with
$B_1(0)=1$. In view of
$(5.9)$ we may terminate the series at a high power with
a negligible error, and we need in fact to deal with
the initial term only, since the oscillating property of other
terms is the same as that of the initial term. The
saddle point method gives that
$$
\eqalignno{
{1\over2\pi\sqrt{y}}&
\int_{-\infty}^\infty \eta_1(t)\left[{t\over2\pi}
\exp\left(-it\log{t\over2\pi ey^{1/r}}\right)\right]^r dt\cr
&=y^{1/2+1/2r}\exp\left(2\pi iry^{1/r}
-\txt{1\over4}\pi i\right)B_2(y),&(5.12)
}
$$
where $B_2(y)$ is a power series of $y^{-1/r}$ 
and the off-diagonal contribution is ignored, as it is
negligible. With this, let  $N$ be large, let
$H=N^\vartheta$, 
with $0<\vartheta<1$, and assume, in place of of $(5.9)$,
that
$$
\hbox{$\log H/\log|\nu_V|$ is large.}\eqno(5.13)
$$
Then, we set, in the above,
$$
w(x)={\phi(x)\over H\sqrt{\pi}}\int_{x-N}^\infty
\exp\big(-(\xi/H)^2\big)d\xi.\eqno(5.14)
$$
Here $\phi(x)$ is a $C^\infty$-function which is equal to
$1$ for $x\in [2H, N+H\log N]$ and has a support in
$(H, N+2H\log N)$. We may
assume that $w^{(q)}(x)\ll H^{-q}$ 
for each fixed $q\ge0$; see [13, Vol.\ II, p.\ 8] for instance.
We have
$$
\sum_{n\le N}\tau_{V\times V}^{(2)}(n)\le
\sum_{n\le N}\tau_{V\times V}^{(2)}(n)w(n)+
O(|\nu_V|^cH),\eqno(5.15)
$$
where the third term is due to $(5.6)$ but with $N$
being replaced by $H$; the non-negativity of 
$\tau_{V\times V}^{(2)}(n)$ plays a r\^ole.
Assertions $(5.7)$--$(5.12)$ imply that the last weighted
sum equals
$$
\eqalignno{
&L(1;\r{sym}^2V)L(1;\r{sym}^4V)\hat{w}(1)\cr
+&\, c\kappa_{V,r}
\sum_{n=1}^\infty{\tau_{V\times V}^{(2)}(n)
\over n^{1/2-1/2r}}\int_0^\infty w(x)x^{1/2r-1/2}
\exp\big(2\pi ir(nx)^{1/r}\big)dx,&(5.16)
}
$$
in which lesser terms have been ignored. We have obviously
$$
\hat{w}(1)=N+O(H\log N).\eqno(5.17)
$$
Applying integration by parts to the last integral many times,
we see that those terms with 
$n\ge (N/H)^r N^{\kappa-1}$ are negligible
for any fixed $\kappa>0$. Thus $(5.16)$ equals
$$
\eqalignno{
&L(1;\r{sym}^2V)L(1;\r{sym}^4V)N+O(|\nu_V|^c H
\log N)\cr+&\, c\kappa_{V,r}
\sum_{n\le (N/H)^r N^{\kappa-1}}
{\tau_{V\times V}^{(2)}(n)
\over n^{1/2-1/2r}}\int_0^\infty w(x)x^{1/2r-1/2}
\exp\big(2\pi ir(nx)^{1/r}\big)dx.&(5.18)
}
$$
When $H=N^{1-1/r+\kappa}$, this sum is empty. 
Our choice $(5.14)$ implies, further, that the last integral
is $\ll N^{1/2-1/2r}n^{-1/r}$. This means, in view of
$(5.6)$, that the last sum is 
$\ll |\nu_V|^c(N/H)^{(r-1)/2}
N^\kappa$. Hence the optimal value 
of $H$ is $N^{(r-1)/(r+1)}$. 
We repeat the same procedure after modifying $(5.14)$
in an obvious way so that $(5.15)$ is replaced 
by the opposite inequality. This ends the proof
of $(3.14)$, since we have $r=9$ actually.
\vskip 1cm
\centerline{\bf References}
\medskip
\noindent
\item{[1]} K. Chandrasekharan and R. Narasimhan.
Functional equations with multiple Gamma factors and
the average order of arithmetical functions. Ann.\ Math.,
{\bf 76} (1962), 93--136.
\item{[2]} S. Gelbart and H. Jacquet. A relation between
automorphic representations of $\r{GL}(2)$ and 
$\r{GL}(3)$. Ann.\ Sci.\ \'Ecole Normale Sup.\ $4^\r{e}$
s\'erie, {\bf 11} (1978), 471--552. 
\item{[3]} S. Gelbart and F. Shahidi. Boundedness
of automorphic $L$-functions in vertical strips.
J. Amer. Math.\ Soc., {\bf 14} (2000), 79--107.
\item{[4]} J. Hoffstein and P. Lockhart. Coefficients of Maass
forms and the Siegel zero. Ann.\ Math., {\bf 140} (1994), 
161--176; D. Goldfeld, J. Hoffstein and D. Lieman.
Appendix: An effective zero-free region. 
ibid, 177--181.
\item{[5]} G. Hoheisel. Primzahlprobleme in der 
Analysis. Sitz.\ Preuss.\ Akad.\ Wiss.\ phys.-math. Klasse, 
1930, pp.\ 580--588.
\item{[6]} H. Kim. Functionality for the exterior square of 
$\r{GL}_4$ and the symmetric fourth 
power of $\r{GL}_2$. J. Amer.\ Math.\
Soc., {\bf 16} (2003), 139--183; Appendix 1 by D. 
Ramakrishnan and Appendix 2 by H. Kim and P. Sarnak.
\item{[7]} H. Kim and F. Shahidi. Cuspidality of symmetric
power with applications. Duke Math., {\bf 112} (2002),
177--197.
\item{[8]} E. Landau. \"Uber die Anzahl der Gitterpunkte in
gewissen Bereichen. Nachr.\ Gesell.\ Wissens.\ G\"ottingen,
math.\ phy. Kl., 1912, Heft 6, 687--770.
\item{[9]} ---: \"Uber die Wurzeln der Zetafunktion. Math. Z., 
{\bf 20} (1924), 98--104.
\item{[10]} H.L. Montgomery. 
{\it Topics in Multiplicative Number 
Theory\/}. Lect.\ Notes in Math., {\bf 227}, 
Springer-Verlag, Berlin 1971.
\item{[11]} Y. Motohashi. Primes in arithmetic progressions.
Invent.\ math., {\bf 44} (1978), 163--178.
\item{[12]} ---. {\it Sieve Methods and 
Prime Number Theory\/}. 
Lect.\ Notes in Math.\ Phys., {\bf72},
Tata IFR and Springer-Verlag, Bombay 1983. 
\item{[13]} ---. {\it Analytic Number Theory\/}.\ Vol.\ I:
{\it Distribution of Prime Numbers\/}. Asakura Publishing, 
Tokyo 2009; Vol.\ II: {\it Zeta Analysis\/}. ibid, 2011. (Japanese)
\item{[14]} ---. Elements of automorphic representations.
arXiv:1112.4226 [math.NT].
\item{[15]} ---. An extension of the Linnik phenomenon.
Proc.\ Steklov Inst.\ Math.,
{\bf 280} (2013), Issue 2 Supplement, 56--64;
Part II. arXiv:1208.6159 [math.NT].
\item{[16]} ---. On sums of Hecke--Maass eigenvalues
squared over primes in arithmetic progressions. In preparation.
\item{[17]} D. Ramakrishnan and S. Wang. 
On the exceptional
zeros of Rankin--Selberg $L$-functions. Compositio Math.,
{\bf 135} (2003), 211--244.
\item{[18]} A. Selberg. Remarks on sieves. In: 
{\it Collected Papers\/}.\ I.
Springer, Berlin 1989, pp.\ 609--615. 
\item{[19]} E.C. Titchmarsh: {\it The Theory of the Riemann
Zeta-Function\/}. Clarendon Press, Oxford 1951.
\vskip 1cm
\noindent 
Honkomagome 5-67-1-901, Tokyo 113-0021, JAPAN
\hfill\def\ymzeta
{\font\brm=cmr17 at 30pt\font\sssrm=cmr5 at 4pt
\font\ssssrm=cmr5 at2.5pt
{{\brm O}\raise 9pt\hbox{\hskip -22pt
$\hfil\raise3pt\hbox{\ssssrm KH}\atop\hbox{
{\sssrm Y}$\zeta$\hskip-1pt{\sssrm M}}$}}}
\ymzeta

\bye